\documentclass[11pt,a4paper]{article}
\usepackage{amsfonts,amssymb,mathtools,amsxtra,amsthm,array,latexsym,multirow,footnote}
\usepackage{geometry}
\newtheorem{definition}{Definition}%[section]
\newtheorem{proposition}[definition]{Proposition}%[section]

\makesavenoteenv{tabular}
\numberwithin{equation}{section}
\numberwithin{definition}{section}
\newcommand{\Ud}{\mathrm{d}}
\newcommand{\Iff}{if\textcompwordmark f}

\DeclareMathOperator{\rank}{rank}
\DeclareMathOperator{\divz}{div}
\DeclareMathOperator{\sign}{sign}

\begin{document}

\title{Ruled surfaces right normalized}

\author{{Stylianos Stamatakis and Ioanna-Iris Papadopoulou}\textbf{\medskip}\\ \emph{Aristotle University of Thessaloniki}\\ \emph{Department of Mathematics}\\ \emph{GR-54124 Thessaloniki, Greece}\\  \emph{e-mail: stamata@math.auth.gr}}
\date{}
\maketitle
\begin{abstract}
\noindent
This paper deals with skew ruled surfaces $\varPhi$ in the Euclidean space $\mathbb{E}^{3}$ which are right normalized, that is they are equipped with relative normalizations, whose support function is of the form $q(u,v) = \frac{f(u) + g(u)\, v}{w(u,v)}$, where $w^2(u,v)$ is the discriminant of the first fundamental form of $\varPhi$. This class of relatively normalized ruled surfaces contains surfaces such that their  relative image $\varPhi^{*}$ is either a curve or it is as well as $\varPhi$ a ruled surface whose generators are, additionally, parallel to those of $\varPhi$. Moreover we investigate various properties concerning the Tchebychev vector field and the support vector field of such ruled surfaces.
\medskip


\end{abstract}

\section{Preliminaries}\label{Sec1}

In this section we present briefly some definitions, results and formulae of relative Differential Geometry of surfaces and Differential Geometry of ruled surfaces in the Euclidean space $\mathbb{E}^{3}$. The reader can use \cite{Pottmann} and \cite{Schirokow} as general references.

In the three-dimensional Euclidean space $\mathbb{E}^{3}$ we denote by $\varPhi=(U,\overline{x})$ a ruled $C^{r}$-surface of nonvanishing Gaussian curvature, $r\geq3$, defined by an injective $C^{r} $-immersion $\overline{x}=\overline{x}(u,v)$ on a region $U:=I\times \mathbb{R}$ ($I\subset \mathbb{R}$  open interval)  of $\mathbb{R}^{2}$.
We denote by $\langle \,,\rangle$ the standard scalar product in $\mathbb{E}^{3}$.
We introduce the so-called standard parameters $u\in I, v \in \mathbb{R}$ of $\varPhi$, such that
    \begin{equation}    \label{1}
        \overline{x}(u,v)=\overline{s}(u) + v\,\overline{e}(u),
    \end{equation}%
and
    \begin{equation}\label{2}
        \left \vert \overline{e}΄\right \vert =|\overline{e}'|=1,\quad \langle \overline{s}',\overline{e}'\rangle =0,
    \end{equation}%
where the differentiation with respect to $u $ is denoted by a prime.
Here $\varGamma: \overline{s}=\overline{s}(u)$ is the striction curve of $\varPhi$ and the parameter $u$ is the arc length along the spherical curve $\overline{e}=\overline{e}(u)$.

The distribution parameter
$$
\delta(u) :=(\overline{s}',\overline{e},\overline{e}'),
$$
the conical curvature
$$
\kappa(u):=(\overline{e},\overline{e}', \overline{e}'')
$$
 and the function
 $$
 \lambda(u) :=\cot \sigma,
 $$
 where
 $$
 \sigma(u) :=\sphericalangle (\overline{e},\overline{s}')
 $$
  is the striction of $\varPhi$,
  $$
  -\frac{\pi }{2}<\sigma \leq \frac{\pi }{2}, \quad \sign \sigma = \sign \delta,
  $$
   are the fundamental invariants of $\varPhi$ and determine uniquely the ruled surface $\varPhi $ up to Euclidean rigid motions.

We also consider the moving frame $\mathcal{D} : = \{\overline{e}, \overline{n}, \overline{z}\}$ of $\varPhi$, where
$$
\overline{n}(u):=\overline{e}'$$
  is the central normal vector  and
  $$
  \overline{z}(u):=\overline{e} \times \overline{n}
  $$
   is the central tangent vector.
It is well known that the following equations
    \begin{equation} \label{10}
    \overline{e}'=\overline{n},\quad \overline{n}'=-\overline{e}+\kappa \,\overline{z},\quad \overline{z}'=-\kappa \,\overline{n}
    \end{equation}
are valid (see \cite[p. 280]{Pottmann}).
Then we have
    \begin{equation}    \label{15}
    \overline{s}'=\delta \,\lambda \, \overline{e}+\delta \, \overline{z}.
    \end{equation}%
We denote partial derivatives of a function (or a vector-valued function) $f$ in the coordinates $u^{1}:=u,\,u^{2}:=v$ by     $f_{/i}, f_{/ij}$ etc.
Then from (\ref{1}) and (\ref{15})\ we obtain
    \begin{equation}\label{20}
    \overline{x}_{/1}=\delta \, \lambda \,\overline{e}+v\,\overline{n}+\delta \,\overline{z},\quad \overline{x}_{/2}=\overline{e},
    \end{equation}
and thus the unit normal vector $\overline{\xi}(u,v) $ to $\varPhi $ is expressed by
    \begin{equation}\label{21}
    \overline{\xi}=\frac{\delta \,\overline{n}-v\,\overline{z}}{w},
    \end{equation}
where $$
w^2 \coloneqq \delta ^{2} + v^{2}
$$
is the discriminant of the first fundamental form of $\varPhi$.
Let $II=h_{ij}\Ud u^i \Ud u^j$ be the second fundamental form of $\varPhi $, where
      \begin{equation}\label{35}
     h_{11}=-\frac{\kappa\, w^{2}+\delta ' \,v-\delta ^{2}\,\lambda }{w},\qquad h_{12}=\frac{\delta }{w}, \qquad  h_{22}=0.
    \end{equation}
The Gaussian curvature $\widetilde{K}(u,v)$ and the mean curvature $\widetilde{H}(u,v)$ of $\varPhi $ are  given by (see \cite{Pottmann})
   \begin{equation}\label{22}
  \widetilde{K}=- \frac{\delta ^{2}}{w^{4}}, \quad \widetilde{H} = -\frac{\kappa w^{2}+\delta ' v+\delta^{2}\lambda}{2w^{3}}.
 \end{equation}
A $C^{s}$-relative normalization of $\varPhi$ is a $C^{s}$-mapping $\overline{y} = \overline{y}(u,v), 1\leq s < r$, defined on $U$, such that
    \begin{equation}\label{45}
    \rank (\{\overline{x}_{/1},\overline{x}_{/2},\overline{y}\})=3,\,\,
    \rank (\{\overline{x}_{/1},\overline{x}_{/2},\overline{y}_{/i}\})=2,\,\,
    i=1,2,\,\,\forall \left(u,v\right) \in U.
    \end{equation}
The pair $\left(\varPhi,\overline {y}\right)$ is called a relatively normalized ruled surface and the line issuing from a point $P \in \varPhi$ in the direction of $\overline{y}$ is called the relative normal of $\varPhi$ at $P$. The pair $\varPhi^{*}=(U,\overline{y})$ is called the relative image of $(\varPhi,\overline{y})$.

The support function of the relative normalization $\bar{y} $  is defined by
$$
q(u,v):=\langle \bar{\xi},\bar{y}\rangle
$$
(see~ \cite{Manhart3}). Because of (\ref{45}), $q $ never vanishes on $U $.
Conversely, when a support function $q$ is given, the relative normalization $\overline{y}$ of the ruled surface $\varPhi$ is uniquely determined and can be expressed in terms of the moving frame $\mathcal{D}$ as follows \cite[p. 179]{Stamatakis3}:
    \begin{equation}\label{70}
    \overline{y}=y_{1}\,\overline{e}+y_{2}\,\overline{n}+y_{3}\,\overline{z},
    \end{equation}
where
    \begin{equation}\label{75}
    y_{1}=-w\frac{\delta q_{/1}+q_{/2}(\kappa \,w^{2}+\delta 'v)}{\delta ^{2}},\quad
    y_{2}=\frac{\delta ^{2}\,q-w^{2}\,v\,q_{/2}}{\delta w},\quad
    y_{3}=-\frac{v\,q+w^{2}\,q_{/2}}{w}.
    \end{equation}
The coefficients $G_{ij}(u,v)$ of the relative metric $G(u,v)$ of $(\varPhi, \overline{y})$, which is indefinite, are given by
$$
G_{ij} = q^{-1} \, h_{ij}.
$$
Then, by taking (\ref{35}) into consideration, the coefficients of the inverse relative metric tensor are computed by
    \begin{equation*}
    G^{(11)}=0,\quad G^{(12)}=\frac{w\,q}{\delta },\quad G^{(22)}= w\,q \, \frac{\kappa\, w^{2}+\delta 'v-\delta ^{2}\,\lambda }{\delta ^{2}}.
    \end{equation*}
For a function (or a vector-valued function) $f$ we denote by $\nabla^{G}\!f$ the first Beltrami differential operator and by $\nabla _{i}^{G}f $ the covariant derivative in the direction $u^i$, both with respect to the relative metric.
The coefficients $A_{ijk}(u,v)$ of the Darboux tensor are defined by
    \begin{equation*}
    A_{ijk}:= q^{-1} \, \langle \overline{\xi},\,\nabla _{k}^{G}\,\nabla _{j}^{G}\,\overline{x}_{/i}\rangle.
    \end{equation*}
Then, by using the relative metric tensor $G_{ij} $ for ``raising and lowering
the indices'', the Pick invariant $J(u,v)$ of $(\varPhi,\overline{y})$ is given by
\begin{equation*}
   J:=\frac{1}{2}A_{ijk}\,A^{ijk}.
\end{equation*}
As we showed in \cite{Stamatakis5} (see equation (2.2)) the Pick invariant is calculated by \begin{equation}\label{80}
    J = \frac{3\left(w^2 q_{/2}+v\,q\right)}{2\delta^2 w^3 \,q}\Big\{w^2 \!\left[ \kappa  q v + 2 \delta q_{/1} +q_{/2} \left(\kappa \,w^2 + \delta ' v-\delta^2 \lambda \right)  \right]  -\delta^2 q \left( \lambda v - \delta '  \right) \Big\}.
\end{equation}
The relative shape operator has the coefficients $B_{i}^{j}(u,v)$ defined by
    \begin{equation}        \label{90}
    \overline{y}_{/i}=:-B_{i}^{j}\, \overline{x}_{/j}.
    \end{equation}
Then, for the relative curvature $K(u,v)$ and the relative mean curvature $H(u,v)$ of $(\varPhi,\overline{y})$ we have
    \begin{equation}    \label{100}
    K:=\det \left(B_{i}^{j}\right),\quad H:=\frac{B_{1}^{1}+B_{2}^{2}}{2}.
    \end{equation}
We mention finally, that among the surfaces $\varPhi \subset \mathbb{E}^3$ with negative Gaussian curvature the ruled surfaces are characterized by the relation
    \begin{equation}\label{136}
    3 H - J -3 S = 0
    \end{equation}
(see \cite{Stamatakis4}), where $S(u,v)$ is the scalar curvature of the relative metric $G$, which is defined formally as the curvature of the pseudo-Riemannian manifold ($\varPhi,G$).
\section{Right normalizations}\label{Sec2}
We focus now our investigation on the main subject of this paper, namely the right normalizations of a skew ruled surface $\varPhi$, that is, relative normalizations which are given by \eqref{70} and \eqref{75} by means of the support function
 \begin{equation}\label{320}
    q=\frac{f+g\,v}{w},
    \end{equation}
where $f$ and $g$ are arbitrary $C^{s+1}$-functions of $u$, such that $q \neq 0$. These normalizations are introduced in \cite{Stamatakis5} by the authors.

When the function $g$ vanishes in $I$, the relative normal at each point $P \in \varPhi$ lies on the corresponding asymptotic plane $\{P;\overline{e},\overline{n}\}$ of $\varPhi$. Normalizations of this type  are called asymptotic and they have been studied by I.~Kaffas and S.~Stamatakis \cite{Stamatakis3}.

Another special case arises when the function $f$ vanishes in $I$. Then the relative normal at each point $P \in \varPhi$ lies on the corresponding central  plane $\{P;\overline{e},\overline{z}\}$ of $\varPhi$. Normalizations of this type are called central and they have been studied   in \cite{Stamatakis5}.

Since both asymptotic and central normalizations belong to the right ones and they have been studied thoroughly in the above mentioned papers, we assume  that in what follows \emph{none of  the functions $f$ and $g$ is vanishing.}

From \eqref{70}, \eqref{75} and \eqref{320} it follows that a right normalization of the given ruled
surface $\varPhi$ is
        \begin{equation}\label{201}
        \overline{y} = \frac{\left(\kappa \, f - \delta \, g' \right) v + \delta' f - \delta \, f' - \delta^2 \, \kappa \, g}{\delta^2} \, \overline{e} + \frac{f}{\delta} \, \overline{n} - g\,\overline{z}.
        \end{equation}
Then, by using \eqref{10}, \eqref{20}, \eqref{90} and \eqref{201}, we obtain the coefficients $B_i^j$ of the relative shape operator of a right normalization:
\begin{align*}
  B_{1}^{1} & = \frac{\delta g'-\kappa f}{\delta^{2}},  \\
    \begin{split}
    B_{1}^{2}& = \frac{1}{\delta^{3}} \Big[\left(2 \kappa \delta' f-\delta \kappa f' - \delta \delta' g'-\delta \kappa' f+ \delta^{2} g'' \right)v+ \delta^{2}f\left( \kappa \lambda +1\right)\\
    & \phantom {a b c d e f} {}{}{}{}{}+ 2\delta '\left(\delta' f-\delta f'\right)
      + \delta^{3}g' \left(\kappa - \lambda\right)+\delta^{3}\kappa'g-\delta \delta'' f+\delta^{2}f'' \Big],
  \end{split}\\
  \quad B_{2}^{1} &= 0,\\
  B_{2}^{2}&= \frac{\delta g'-\kappa f}{\delta^{2}}.
\end{align*}
Hence, via  \eqref{100}, the relative mean curvature $H$ and the relative curvature $K$  are
        \begin{equation}\label{202}
            H = \frac {\delta \, g' - \kappa \, f} {\delta^{2}}, \quad K = H^{2}.
        \end{equation}
Firstly, we observe that \emph{all points of $\varPhi $ are relative umbilics} ($H^{2}-K\equiv0 $). Thus, for the relative principal curvatures $k_{1}$ and $k_{2}$, which by definition are the eigenvalues of the relative shape operator (see \cite[p. 215]{Schirokow}),
$$
k_{1} = k_{2} = H
$$
holds.

Then, from \eqref{80} we find for the Pick invariant
 \begin{equation}\label{210}
 J = 3g \, \frac {\kappa\, g\, v^{2}+ 2 \delta \,g 'v+\delta^{2} g \,(\kappa-\lambda)-\delta'f+ 2\delta f'} {2\delta^{2} \left(f + g \, v\right)}.
 \end{equation}
Consequently $J$ vanishes identically \Iff{}
\begin{equation*}
\kappa\, g\, v^{2}+ 2 \delta \,g 'v+\delta^{2} g \,(\kappa-\lambda)-\delta'f+ 2\delta f' = 0,
\end{equation*}
or, equivalently, after successive differentiations of this last equation relative to $v$, \Iff{}
\begin{equation}\label{215}
\kappa = g ' = \delta^{2} g \,(\kappa-\lambda)-\delta'f+ 2\delta f' = 0,
\end{equation}
from which we have $$\kappa = 0,$$ i.e., $\varPhi $ is conoidal, $$g = c_{1} \in \mathbb{R}^{*}$$ and $$f = |\delta|^{1/2} \left(\! \frac{c_{1}}{2} \!\! \int \! |\delta|^{1/2} \, \lambda \, \Ud u + c_{2} \right),\, c_{2}\in\mathbb{R}.$$
Thus, the following has been shown
\begin{proposition}
The Pick invariant of a right normalized skew ruled surface $\varPhi \subset \mathbb{E}^{3}$ vanishes identically \Iff{} $\varPhi $ is conoidal, the function $g$ is a nonvanishing constant $c_{1}$ and the function $f$ is given by
 \begin{equation*}
 f = |\delta|^{1/2} \left(\! \frac{c_{1}}{2} \!\! \int \! |\delta|^{1/2} \, \lambda \, \Ud u + c_{2} \right), \, c_{2}\in\mathbb{R}.
\end{equation*}
\end{proposition}
Additionally, in view of (\ref{202}a) and \eqref{215}, \emph{a right normalized ruled surface with vanishing Pick invariant is relatively minimal.}

By using \eqref{136}, (\ref{202}a) and \eqref{210} we obtain the scalar curvature of the relative metric
        \begin{equation*}
        S = -\frac { \kappa\, g ^{2}v^{2}+2 \kappa\, f\, g\, v + \delta^{2} g^{2} (\kappa-\lambda) + 2\kappa\, f^{2}-\delta' f \,g + 2\delta\, (f'g-f\,g')} {2 \delta^{2} \left( f + g \, v \right)}.
        \end{equation*}
The scalar curvature of the relative metric $G$ vanishes identically \Iff{}
\begin{equation*}
  \kappa =\delta^{2} g^{2} (\kappa-\lambda) + 2\kappa\, f^{2}-\delta' f \,g + 2\delta\, (f'g-f\,g')=0,
\end{equation*}
that is, \Iff{} $$\kappa=0$$ and $$ f =\frac {1}{2} |\delta|^{1/2} \, g \left( \int \! |\delta|^{1/2} \, \lambda \, \Ud u +c\right), \, c\in\mathbb{R}.$$
So, we have:
\begin{proposition}\label{A}
The scalar curvature $S$ of a right normalized skew ruled surface $\varPhi \subset \mathbb{E}^{3}$ vanishes identically \Iff{} $\varPhi $ is conoidal and the function $f$ is given by
\begin{equation*}
         f = \frac {1}{2}|\delta|^{1/2} \, g \left( \int \! |\delta|^{1/2} \, \lambda \, \Ud u +c\right), \, c\in\mathbb{R}.
  \end{equation*}
\end{proposition}
We distinguish the right normalizations in two types.
\subsection{Right normalizations of the first type}\label{SubSec21}
We say that a right relative normalization $\overline {y}$ is of the first type if the relative image $\varPhi^*$  of $\left(\varPhi,\overline {y}\right)$ degenerates into a curve. Obviously this occurs \Iff{}
$$
\delta \, g' - \kappa \, f = 0
$$ (cf. \eqref{201}). Thus, on account of \eqref{201} and (\ref{202}a), we conclude:
\begin{proposition}
Let $(\varPhi,\overline{y})$ be a right normalized ruled surface. Then
the following properties are equivalent:\\
\textup{(a)} $\overline{y}$ is a right normalization of the first type.\\
\textup{(b)} $(\varPhi,\overline{y})$ is relatively minimal.\\
\textup{(c)} The function $g$ is given by
        \begin{equation*}
        g = \int \!\frac {\kappa \, f} {\delta} \, \Ud u +c, \, c\in  \mathbb{R}.
        \end{equation*}
\end{proposition}
The right normalized ruled surfaces with vanishing Pick invariant belong obviously to this subclass.

The relative image $\varPhi^*$ is the curve parametrized by
\begin{equation*}
  \overline{y} = \frac{\delta' f - \delta  f' - \delta^2 \, \kappa \, g}{\delta^2} \, \overline{e} + \frac{f}{\delta} \, \overline{n} - g\,\overline{z}.
\end{equation*}
\subsection{Right normalizations of the second type}\label{SubSec22}
A right relative normalization $\overline {y}$ is said to be  of the second  type if the relative image $\varPhi^*$  of $\left(\varPhi,\overline {y}\right)$ does not degenerate into a curve of $\mathbb{E}^3$.
Then $\varPhi^{*}$ is %as well as $\varPhi$
a ruled surface whose generators are parallel to those of $\varPhi$.
%Consequently $\varPhi$  and $\varPhi^{*}$ posses common conical curvature.
From \eqref{201} we find the following parametrization of the striction curve of $\varPhi^*$:
\begin{equation*}
  \varGamma^{*} \colon \overline{s}^{*}=\frac{\delta'  f - \delta f' - \delta^{2} \kappa \, g}{\delta^{2}}\,\overline{e}+\frac{f}{\delta}\,\overline{n}-g\, \overline{z}.
\end{equation*}
Consequently $\varPhi^{*}$ can be parametrized like \eqref{1} and \eqref{2}:
\begin{equation*}
  \varPhi^{*} \colon \overline{y}= \overline{s}^{*} + v^{*}\, \overline{e},
\end{equation*}
where $$v^{*} \coloneqq \frac{(\kappa f - \delta \,g')v}{\delta^{2}}.$$
Considering $\mathcal{D}$ as moving frame of $\varPhi^{*}$  we compute its fundamental invariants:
\begin{equation*}
  \kappa^{*}=\kappa, \, \delta^{*}=\frac{\kappa f-\delta g'}{\delta}, \,  \lambda^{*}=-\frac{\delta^{3}(\kappa g'+\kappa' g)+\delta^{2}(f+f'')-\delta(\delta''f+2 \delta'f')+2\delta'^{2}f}{\delta^{2}(\kappa f-\delta g')}.
\end{equation*}
By using (\ref{21}b) we infer that
\begin{equation*}
  w^{*}=|H|\, w
\end{equation*}
and, thus, by means of (\ref{22}a), the Gaussian curvature $\widetilde{K}^{*}$ of $\varPhi^{*}$ is
\begin{equation*}
  \widetilde{K}^{*}=-\frac{\delta^{6}}{w^{4}(\kappa f-\delta \, g')^{2}}.
\end{equation*}
The focal surfaces, which are the loci of the edges of regression of the developable surfaces consisting of the relative normals along the relative lines of curvature, coincide. The parametrization of the unique relative focal surface of $\varPhi$, which initially reads
        \begin{equation*}
        \overline{x}^{*}=\overline{s}+v \, \overline{e}+\frac{1}{H} \,\overline{y},
        \end{equation*}
in view of \eqref{201} and (\ref{202}a) becomes
        \begin{equation*}
        \overline{x}^{* } = \overline{s} + \frac {\left( \delta'  f - \delta  f' - \delta^{2} \, \kappa \, g \right) \overline{e} + \delta \, f \, \overline{n} - \delta^{2}  g \, \overline{z}} {\delta \, g' - \kappa \, f},
        %\overline{x}^{* } = \overline{s}+ \frac{1}{H} \, \overline{s}^*
        \end{equation*}
i.e., \emph{the focal surface degenerates into a curve $\varLambda ^{* } $ and all relative normals along each generator form a pencil of straight lines}.

\section{The Tchebychev vector field of a right normalization}\label{Sec3}
In \cite{Stamatakis3} it was shown that the coordinate functions of the Tchebychev vector $\overline{T}(u,v)$ of $(\varPhi,\overline{y})$, which is defined by
    \begin{equation*}
    \overline{T}:=T^{m}\, \overline{x}_{/m},
    \end{equation*}
where $$T^{m}:=\frac{1}{2}A_{i}^{im},$$
are given by
    \begin{equation}\label{121}
    T^{1}=\frac{w^{2}q_{/2}+vq}{\delta w},\,\, T^{2}=\frac{2\delta\,w^{2}q_{/1}+\delta'q(\delta^{2}-v^{2})}{2\delta^{2}\, w}+\frac{T^{1}(\kappa w^{2}+\delta'v-\delta^{2}\lambda)}{\delta}.
    \end{equation}
By means of \eqref{20} and \eqref{320} the Tchebychev vector of a right normalization can be expressed in terms of the moving frame $\mathcal{D}$ as follows:
\begin{equation}\label{235}
  \overline{T}=\frac{2\kappa\, g\, v^{2}+ (\delta'g+2\delta\, g')v + 2\delta^{2}\kappa\,g - \delta'f+2\delta f'}{2\delta^{2}}\,\overline{e}+\frac{g}{\delta}(v\,\overline{n}+\delta\,\overline{z}).
\end{equation}
The vectors $\overline{T}$ are orthogonal to the generators \Iff{} $\langle \overline{e},\overline{T}\rangle=0$. Taking \eqref{235} into consideration we find
\begin{equation*}
   2 \kappa \,g\, v^{2}+(\delta'g+2\delta \,g')v + 2\delta^{2}\kappa\,g -\delta'f+2\delta f'=0,
\end{equation*}
or, after successive differentiations of this last equation relative to $v$, \Iff{}
\begin{equation*}
    2 \kappa \,g =\delta'g+2\delta \,g' = 2\delta^{2}\kappa\,g -\delta'f+2\delta f' = 0.
\end{equation*}
After standard treatment of this system we deduce that $$\kappa=0,$$ $$g = c_{1}|\delta|^{-1/2}, \,c_{1}\in\mathbb{R^{*}}$$ and $$f= c_{2}|\delta|^{1/2}, \, c_{2}\in\mathbb{R^{*}}. $$ So, we have the following
 \begin{proposition}
   The Tchebychev vector field $ \overline{T}$ of a right normalized skew ruled surface $\varPhi \subset \mathbb{E}^{3}$ is orthogonal to the generators of $\varPhi$ \Iff{} $\varPhi$ is conoidal and the functions $g$ and $f$ are given by
   \begin{equation*}
     g= c_{1}|\delta|^{-1/2}, c_{1}\in\mathbb{R^{*}}\,\, \textit{and}\,\, f= c_{2}|\delta|^{1/2}, c_{2}\in\mathbb{R^{*}}.
   \end{equation*}
 \end{proposition}
We turn now to the right normalized ruled surfaces $\left(\varPhi,\overline {y}\right)$, whose Tchebychev vectors are tangent or orthogonal to one of the following geometrically distinguished families of curves of $\varPhi$:
\begin{description}
  \item[a.] the curves of constant striction distance ($u$-curves),
  \item[b.] the curved asymptotic lines and
  \item[c.] the $\widetilde{K}$-curves, i.e., the curves along which the Gaussian curvature is constant \cite{Sachs}.
\end{description}
The corresponding differential equations of these families of curves are
\begin{align}
  &v'=0, \label{545} \\
  &\kappa \,v^{2}+\delta'v+\delta^{2}(\kappa-\lambda)-2\delta \,v'= 0, \label{550} \\
  & 2 \delta \,v\, v'+ \delta'\left(\delta^{2}-v^{2}\right)=0.\label{555}
\end{align}
We will investigate necessary and sufficient conditions for the Tchebychev vector field $\overline{T}$ to be tangential or orthogonal to each one of these families of curves.

We consider a directrix $ \varLambda : v=v(u)$ of $\varPhi$. Then we have
\begin{equation}\label{560}
  \overline{x}'=(\delta\,\lambda+v')\,\overline{e}+v\, \overline{n}+\delta\,\overline{z}.
\end{equation}
From \eqref{235} and \eqref{560} it follows: $ \overline{x}'$ and $\overline{T}$ are parallel or orthogonal \Iff{}
\begin{equation}\label{565}
2 \kappa \,g\, v^{2} + \left(\delta' g+ 2 \delta \,g'\right)v + 2\delta^{2}\kappa\,g -\delta'f+ 2\delta f'-2 \delta g \left(\delta\,\lambda+v'\right)=0
\end{equation}
or
\begin{equation}\label{570}
  \Big[2 \kappa g  v^{2} + \left(\delta' g+ 2 \delta g'\right)v+ 2\delta^{2}\kappa\,g -\delta'f+ 2\delta f'\Big]\left(\delta\,\lambda+v'\right)+ 2\delta\, g\, w^{2}=0,
\end{equation}
respectively.

From \eqref{545} and \eqref{565}, resp. \eqref{570}, we have: $\overline{T}$ is tangential or orthogonal to
 the $u$-curves \Iff{}
 \begin{equation}\label{549}
   2 \kappa\, g\, v^{2}+\left(\delta' g + 2 \delta g'\right)v+ 2 \delta^{2}g(\kappa-\lambda)-\delta' f + 2 \delta f'=0
 \end{equation}
 or
 \begin{equation}\label{551}
   2g (\kappa\,\lambda+1)v^{2}+\lambda (\delta' g + 2 \delta g')v + 2\delta^{2}g (\kappa\,\lambda+1)+\lambda(2 \delta f'-\delta' f)=0,
 \end{equation}
 respectively.
 From \eqref{549} we find that $\overline{T}$ is tangential to the $u$-curves \Iff{}
 \begin{equation*}
   \kappa=\delta' g + 2 \delta g'=2 \delta^{2}g(\kappa-\lambda)-\delta' f + 2 \delta f'=0,
 \end{equation*}
 that is, \Iff{} $$\kappa=0,$$ $$g= c_{1}|\delta|^{-1/2},\, c_{1}\in\mathbb{R^{*}}$$ and $$f=|\delta|^{1/2}\left(c_{1}\int\lambda \, \Ud u+c_{2}\right),\, c_{2}\in\mathbb{R}.$$

 From \eqref{551} we derive that $\overline{T}$ is orthogonal to the $u$-curves \Iff{}
  \begin{equation*}
    \kappa\,\lambda+1=\lambda (\delta' g + 2 \delta\, g')=2\delta^{2}g (\kappa\,\lambda+1)+\lambda(2 \delta \,f'-\delta' f)=0.
  \end{equation*}
  By direct computation we deduce that $$\kappa \,\lambda+1=0,$$ i.e., the striction curve of $\varPhi$ is an Euclidean line of curvature, $$g= c_{1}|\delta|^{-1/2},\, c_{1}\in\mathbb{R^{*}}$$ and $$f= c_{2}|\delta|^{1/2},\, c_{2}\in\mathbb{R^{*}}.$$
 Therefore, we obtain
   \begin{proposition}
   The Tchebychev vector field $ \overline{T}$ of a right normalized skew ruled surface $\varPhi \subset \mathbb{E}^{3}$ is\\
   \textup{(a)} tangential to the $u$-curves of $\varPhi$ \Iff{} $\varPhi$ is conoidal and the functions $g$ and $f$ are given by
     \begin{equation*}
       g= c_{1}|\delta|^{-1/2}, \,c_{1}\in\mathbb{R^{*}}\, \, \textit{and}\,\, f=|\delta|^{1/2}\left(c_{1}\int\lambda \, \Ud u+c_{2}\right),\,c_{2}\in\mathbb{R}.\\
     \end{equation*}
     \textup{(b)}orthogonal to the $u$-curves of $\varPhi$ \Iff{} the striction curve of $\varPhi$ is an Euclidean line of curvature and the functions $g$ and $f$ are given by
     \begin{equation*}
       g= c_{1}|\delta|^{-1/2},\, c_{1}\in\mathbb{R^{*}}\,\,\textit{and}\,\,f= c_{2}|\delta|^{1/2},\,c_{2}\in\mathbb{R^{*}}.
     \end{equation*}
     \end{proposition}
From \eqref{545} and \eqref{565} we infer, that $\overline{T}$ is tangential to the curved asymptotic lines \Iff{}
\begin{equation*}
  \kappa g v^{2}+ 2 \delta g' v + \delta^{2} g (\kappa-\lambda)- \delta'f+ 2\delta f'=0,
\end{equation*}
that is, \Iff{}
\begin{equation*}
  \kappa = g'=\delta^{2} g (\kappa-\lambda)- \delta'f+ 2\delta f'=0,
\end{equation*}
from which we have $$\kappa=0,$$ $$g = c_{1} \in \mathbb{R}^{*}$$ and $$f=|\delta|^{1/2}\left(\frac{c_{1}}{2}\int|\delta|^{1/2}\lambda \,\Ud u+c_{2}\right),\,c_{2}\in\mathbb{R}.$$
So, we arrive at
  \begin{proposition}
   The Tchebychev vector field $ \overline{T}$ of a right normalized skew ruled surface $\varPhi \subset \mathbb{E}^{3}$ is tangential to the curved asymptotic lines of $\varPhi$ \Iff{} $\varPhi$ is conoidal, the function $g$ is a nonvanishing constant $c_{1}$ and the function $f$ is given by
   \begin{equation*}
   f=|\delta|^{1/2}\left(\frac{c_{1}}{2}\int|\delta|^{1/2}\lambda \, \Ud u+c_{2}\right),\,c_{2}\in\mathbb{R}.
   \end{equation*}
 \end{proposition}
 From \eqref{555} and \eqref{565}, resp. \eqref{570}, we infer: $\overline{T}$ is tangential or orthogonal to
 the $\widetilde{K}$-curves \Iff{}
 \begin{equation}\label{552}
   2 \kappa g v^{3}+ 2 \delta g' v^{2}+\left[2 \delta^{2}g(\kappa-\lambda)-\delta' f+ 2 \delta f'\right]v+\delta^{2}\delta' g=0
 \end{equation}
 or
 \begin{equation}\label{553}
   \begin{split}
    & 2 \kappa \delta' g v^{4}+ \left[ 4\delta^{2}g(\kappa\lambda+1)+\delta'(\delta' g+ 2 \delta g')\right]v^{3} \\
    &+\left(2 \delta^{2}\delta'\lambda g-\delta'^{2}f+ 2 \delta \delta'f'+4\delta^{3}\lambda g'\right)v^{2}\\
    &+\delta^{2}\left[4\delta^{2}g(\kappa\lambda+1)-2\delta'\lambda f-\delta'^{2}g+4\delta\lambda f'-2\delta \delta'g'\right]v\\
    &-\delta^{2}\delta'\left(2 \delta^{2}\kappa g-\delta'f+2 \delta f'\right)=0,
   \end{split}
 \end{equation}
respectively.
From \eqref{552} we find that $\overline{T}$ is tangential to the $\widetilde{K}$-curves \Iff{}
\begin{equation*}
  \kappa = g'=2 \delta^{2}g(\kappa-\lambda)-\delta' f+ 2 \delta f'=\delta'=0,
\end{equation*}
i.e., \Iff{} $$\kappa=0,$$ $$\delta=c_{1}\in \mathbb{R}^{*},$$ $$g = c_{2} \in \mathbb{R}^{*}$$ and $$f= c_{1} c_{2} \int\lambda \, \Ud u+c_{3},\, c_{3}\in\mathbb{R}.$$ From \eqref{553} we deduce that $\overline{T}$ is orthogonal to the $\widetilde{K}$-curves \Iff{}
\begin{equation*}
\begin{split}
    & \kappa \delta'=4\delta^{2}g(\kappa\lambda+1)+\delta'(\delta' g+ 2 \delta g')=2 \delta^{2}\delta'\lambda g-\delta'^{2}f+ 2 \delta \delta'f'+4\delta^{3}\lambda g'=0, \\
     &4\delta^{2}g(\kappa\lambda+1)-2\delta'\lambda f-\delta'^{2}g+4\delta\lambda f'-2\delta \delta'g'=\delta'\left(2 \delta^{2}\kappa g-\delta'f+2 \delta f'\right)=0,
\end{split}
\end{equation*}
that is, \Iff{} $$\delta=c\in \mathbb{R}^{*}$$ or $$\kappa=0.$$ If $$\delta=c\in \mathbb{R}^{*},$$ we deduce that $$\kappa\lambda+1=0,$$ i.e., $\varPhi$ is an Edlinger surface \footnote{i.e., a ruled surface whose osculating quadrics are rotational hyperboloids. The Edlinger surfaces are characterized by the conditions $\delta' = \kappa  \, \lambda +1 = 0$ (see \cite[p. 36]{Hoschek}, \cite{Sachs}).}, $$g = c_{1} \in \mathbb{R}^{*}$$ and $$f = c_{2} \in \mathbb{R}^{*}.$$ If $$\kappa=0$$ and $$\delta \neq c\in \mathbb{R}^{*}$$ we arrive at a contradiction.
Thus, the following has been shown
 \begin{proposition}
   The Tchebychev vector field $ \overline{T}$ of a right normalized skew ruled surface $\varPhi \subset \mathbb{E}^{3}$ is\\
   \textup{(a)} tangential to the $\widetilde{K}$-curves of $\varPhi$ \Iff{} $\varPhi$ is conoidal of constant distribution parameter $c_{1}$, the function $g$ is a nonvanishing constant $c_{2}$ and the function f is given by
     \begin{equation*}
      f= c_{1} c_{2} \int\lambda \, \Ud u+c_{3}, \, c_{3}\in\mathbb{R}.\\
     \end{equation*}
   \textup{(b)}  orthogonal to the $\widetilde{K}$-curves of $\varPhi$ \Iff{} $\varPhi$ is an Edlinger surface and the functions $g$ and $f$ are nonvanishing constants $c_{1}$ and $c_{2}$, respectively.
    \end{proposition}
The following table summarizes the results:
\begin{footnotesize}
\begin{center}
  \begin{tabular}{ | m{3cm} | m{2.7cm} | m{2cm} | m{4.4cm} | }
    \hline
     \centering $ \overline{T}$ is \dots & \centering Type of the ruled surface $\varPhi$   & \centering $g$ & \centering \arraybackslash $f$ \\ [0.1cm]
     \hline \hline
    orthogonal to the generators & conoidal & $g= c_{1}|\delta|^{-1/2}$,

    $ c_{1}\in\mathbb{R^{*}}$ & $f= c_{2}|\delta|^{1/2}, \,c_{2}\in\mathbb{R^{*}} $ \\
    \hline
     tangential to the

    $u$-curves & conoidal & $g= c_{1}|\delta|^{-1/2}$,

     $c_{1}\in\mathbb{R^{*}}$ & $f=|\delta|^{1/2}\left(c_{1}\int\lambda \, \Ud u+c_{2}\right)$,

    $ c_{2}\in\mathbb{R} $ \\
    \hline
    orthogonal to the

    $u$-curves & the striction curve is an Euclidean line of curvature & $g= c_{1}|\delta|^{-1/2}$,

    $c_{1}\in\mathbb{R^{*}}$ & $f= c_{2}|\delta|^{1/2}, \, c_{2}\in\mathbb{R^{*}} $ \\
    \hline
    tangential to the curved asympt. lines & conoidal & $g \!=\! c_{1} \!\in \! \mathbb{R}^{*}$ & $f\!=\!|\delta|^{1/2}\!\left(\frac{c_{1}}{2}\!\int \!|\delta|^{1/2}\lambda \,\Ud u+c_{2}\right)$,

    $c_{2}\!\in \!\mathbb{R} $ \\
    \hline
    tangential to the

    $\widetilde{K}$-curves & conoidal,

    $\delta=c_{1}\in \mathbb{R}^{*}$  & $g = c_{2} \in \mathbb{R}^{*}$ & $f= c_{1} c_{2} \int\lambda \, \Ud u+c_{3},\, c_{3}\in\mathbb{R} $ \\
    \hline
    orthogonal to the

    $\widetilde{K}$-curves & Edlinger surface & $g = c_{1} \in \mathbb{R}^{*}$ & $f = c_{2} \in \mathbb{R}^{*}$ \\
    \hline
  \end{tabular}
\end{center}
\end{footnotesize}
The divergence $\divz^{I}\overline{T}$ of $\overline{T}$ with respect to the first fundamental form $I$ of $\varPhi$, which initially reads (see \cite{Stamatakis3})
\begin{equation*}
\divz^{I}\overline{T}=\frac{\left(w T^{i} \right)_{/i}}{w}
\end{equation*}
becomes, on account of \eqref{121} and \eqref{320},
\begin{equation*}
  \divz^{I}\overline{T}=\frac{6 \kappa g v^{3}+6 \delta g' v^{2}+\left(6 \delta^{2}\kappa g-2\delta^{2}\lambda g-\delta'f+2\delta \,f'\right)v+\delta^{2}\left(\delta'g+4 \delta \,g'\right)}{2 \delta^{2}w^{2}},
\end{equation*}
from which we have that the Tchebychev vector field $ \overline{T}$ is incompressible with respect to the first fundamental form of $\varPhi$ $(\divz^{I}\overline{T}=0)$ \Iff{}
\begin{equation*}
 \kappa= g'=6 \delta^{2}\kappa g-2\delta^{2}\lambda g-\delta'f+2\delta f'=\delta'g+4 \delta g'=0,
\end{equation*}
or \Iff{} $$\kappa=0,$$ $$g=c_{1}\in \mathbb{R^{*}}, $$ $$\delta=c_{2}\in \mathbb{R^{*}}$$  and $$ f=c_{1} c_{2} \int \lambda \, \Ud u + c_{3},\, c_{3}\in \mathbb{R}.$$ Therefore, we arrive at
 \begin{proposition}
   The Tchebychev vector field $ \overline{T}$ of a right normalized skew ruled surface $\varPhi \subset \mathbb{E}^{3}$ is incompressible with respect to the first fundamental form of $\varPhi$ \Iff{} $\varPhi$ is conoidal of constant distribution parameter $c_{2}$, the function $g$ is a nonvanishing constant $c_{1}$ and the function $f$ is given by
    \begin{equation*}
      f=c_{1} c_{2} \int \lambda \, \Ud u + c_{3},\, c_{3}\in \mathbb{R}.
    \end{equation*}
    \end{proposition}
Let, now, $\divz^{G}\overline{T}$ be the divergence of $\overline{T}$ with respect to the relative metric of $(\varPhi, \overline{y})$. Analogously to the above computation, by using \eqref{35}, we get
\begin{equation*}
  \divz^{G}\overline{T}=\frac{\kappa g^{2}v^{2}+2\kappa f g v- \delta^{2}g^{2}\left(\kappa-\lambda\right)+\delta'f g-2\delta g f'+2\delta f g'}{\delta^{2}\left(gv+f\right)}.
\end{equation*}
The Tchebychev vector field $ \overline{T}$ is incompressible with respect to the relative metric $(\divz^{G}\overline{T}=0)$ \Iff{}
\begin{equation*}
  \kappa=- \delta^{2}g^{2}\left(\kappa-\lambda\right)+\delta'f g-2\delta g f'+2\delta f g'=0,
\end{equation*}
i.e., \Iff{} $$\kappa=0$$ and $$f=\frac{1}{2}|\delta|^{1/2}g\left(\int|\delta|^{1/2} \lambda \, \Ud u+c \right),\, c\in \mathbb{R}.$$
So, by taking into consideration Proposition \ref{A}, we deduce:
 \begin{proposition}
 Let $\varPhi \subset \mathbb{E}^{3}$ be a right normalized skew ruled surface. The following properties are equivalent:\\
 \textup {(a)} The Tchebychev vector field $ \overline{T}$ is incompressible  with respect to the relative metric.\\
 \textup {(b)} The scalar curvature $S$ of the relative metric vanishes identically.\\
 \textup {(c)} $\varPhi$ is conoidal and the function $f$ is given by
   \begin{equation*}
     f=\frac{1}{2}|\delta|^{1/2}g\left(\int|\delta|^{1/2} \lambda \, \Ud u+c\right),\, c\in \mathbb{R}.
   \end{equation*}
    \end{proposition}
\section{The support vector field of a right normalization}\label{Sec4}
Let
\begin{equation*}
  \overline{Q}:= \frac{1}{4} \bigtriangledown^{G}\!\Big(\frac{1}{q},\overline{x}\Big)
\end{equation*}
be the support vector $\overline{Q}(u,v)$ of $( \varPhi, \overline{y})$, which is introduced in  \cite{Stamatakis3}.
On account of \eqref{20},  \eqref{35} and \eqref{320} we express the support vector in terms of the moving frame $\mathcal{D}$ as follows:
\begin{equation}            \label{571}
\overline{Q}=-w \frac{\left(\delta\, g'-\kappa\, f\right)v+\delta^{2}\kappa \,g- \delta ' f+\delta f'}{4\delta^{2}(g\,v+f)}\overline{e}+\frac{f\,v-\delta^{2}g}{4\delta \,w(g\,v+f)}\left(v\, \overline{n}+\delta\,\overline{z}\right).
\end{equation}
 The vectors $\overline{Q}$ are orthogonal to the generators \Iff{} $\langle \overline{e},\overline{Q}\rangle=0$. Taking \eqref{571} into consideration we have
\begin{equation*}
  \left(\delta g'-\kappa f\right)v+\delta^{2}\kappa g -\delta'f+\delta f' =0,
\end{equation*}
that is, \Iff{}
\begin{equation*}
  \delta g'-\kappa f=\delta^{2}\kappa g -\delta'f+\delta f'=0,
\end{equation*}
from which we find that $\varPhi$ is relative minimal and $$f=\pm\delta\left|c-g^{2}\right|^{1/2},\, c\in \mathbb{R},\, g^{2}\neq c.$$
Thus, we arrive at:
 \begin{proposition}
   The support vector field $\overline{Q}$ of a right normalized skew ruled surface $\varPhi \subset \mathbb{E}^{3}$ is orthogonal to the generators of $\varPhi$ \Iff{} $\varPhi$ is relative minimal and the function $f$ is given by
   \begin{equation*}
     f=\pm\delta\left|c-g^{2}\right|^{1/2},\, c\in \mathbb{R},\, g^{2}\neq c.
   \end{equation*}
   \end{proposition}
We will investigate, now, the right normalized ruled surfaces $\varPhi$, whose support vectors are tangent or orthogonal to the above mentioned geometrically distinguished families of curves of $\varPhi$.
From \eqref{560} and \eqref{571} we have: $\overline{x}'$ and $\overline{Q}$ are parallel or orthogonal \Iff{}
\begin{equation}\label{572}
w^{2}\left[\left(\delta g' -  \kappa f\right)v+\delta^{2}\kappa g - \delta' f+ \delta f'\right]+\delta\left(fv-\delta^{2}g\right)\left(\delta \lambda + v'\right)=0
\end{equation}
or
\begin{equation}\label{573}
-\left(\delta \lambda + v'\right)\left[\left(\delta g' -  \kappa f\right)v+\delta^{2}\kappa g - \delta' f+ \delta f'\right]+\delta\left(fv-\delta^{2}g\right)=0.
\end{equation}
From \eqref{550} and \eqref{572}, resp. \eqref{573}, we find: $\overline{Q}$ is tangential or orthogonal to the $u$-curves \Iff{}
\begin{equation}\label{580}
  \begin{split}
      & \left(\kappa f -  \delta g'\right)v^{3}+\left(-\delta^{2} \kappa g+ \delta' f- \delta f'\right)v^{2}+ \delta^{2} \left[f \left(\kappa-\lambda\right) -\delta g'\right]v\\
       & -\delta^{2}\left[\delta^{2}g \left(\kappa-\lambda\right)- \delta'f+ \delta f'\right]=0
  \end{split}
\end{equation}
or
\begin{equation}\label{585}
  \left[f\left(\kappa \lambda + 1\right)-\delta\lambda g'\right]v-\delta^{2}g \left(\kappa\lambda+1\right)+\lambda \left(\delta'f-\delta f' \right)=0,
\end{equation}
respectively. From \eqref{580} we infer that $\overline{Q}$ is tangential to the $u$-curves \Iff{}
\begin{equation*}
  \kappa f -  \delta g'=-\delta^{2} \kappa g+ \delta' f- \delta f'=f \left(\kappa-\lambda\right) -\delta g'=\delta^{2}g \left(\kappa-\lambda\right)- \delta'f+ \delta f'=0,
\end{equation*}
that is, \Iff{} $\varPhi$ is relative minimal, $$\lambda=0,$$ i.e., $\varPhi$ is orthoid \footnote{that is, a ruled surface whose striction curve is an orthogonal trajectory of the generators. The ortoid ruled surfaces are characterized by the condition $  \lambda  = 0$.}  and $$f=\pm \delta \left|c-g^{2}\right|^{1/2},\, c \in \mathbb{R},\,g^{2}\neq c.$$ From \eqref{585} we take that $\overline{Q}$ is orthogonal to the $u$-curves \Iff{}
\begin{equation*}
  f\left(\kappa \lambda + 1\right)-\delta\lambda g'= -\delta^{2}g \left(\kappa\lambda+1\right)+\lambda \left(\delta'f-\delta f' \right)=0,
\end{equation*}
i.e., \Iff{} $$\kappa \lambda +1=\frac{\delta \lambda  g'}{f}$$ and $$f=\pm \delta \left|c-g^{2}\right|^{1/2},\, c \in \mathbb{R},\, g^{2}\neq c,$$ hence $$\kappa= \pm g' \left|c-g^{2}\right|^{-1/2}- \lambda^{-1},\, \lambda\neq0.$$
Therefore, we obtain
 \begin{proposition}
   The support vector field $\overline{Q}$ of a right normalized skew ruled surface $\varPhi \subset \mathbb{E}^{3}$ is\\
   \textup{(a)}tangential to the $u$-curves of $\varPhi$ \Iff{} $\varPhi$ is an orthoid, relative minimal surface and the function $f$ is given by
         \begin{equation*}
         f=\pm \delta \left|c-g^{2}\right|^{1/2},\, c \in \mathbb{R},\, g^{2}\neq c.\\
         \end{equation*}
   \textup{(b)}orthogonal to the $u$-curves of $\varPhi$ \Iff{} the conical curvature and the function $f$ are given by
     \begin{equation*}
      \kappa= \pm g' \left|c-g^{2}\right|^{-1/2}- \lambda^{-1},\,c \in \mathbb{R},\, \lambda\neq0,\, g^{2}\neq c \,\, \textit{and}\,\,f=\pm \delta \left|c-g^{2}\right|^{1/2}.
     \end{equation*}
  \end{proposition}
From \eqref{545} and \eqref{572} we have, that $\overline{Q}$ is tangential to the curved asymptotic lines \Iff{}
\begin{equation*}
\begin{split}
   & \left(\kappa f - 2 \delta g'\right)v^{3}+\left(-\delta^{2} \kappa g+ \delta' f- 2\delta f'\right)v^{2}+ \delta^{2} \left[f \left(\kappa-\lambda\right)+\delta' g -2\delta g'\right]v \\
     & -\delta^{2}\left[\delta^{2}g \left(\kappa-\lambda\right)-2 \delta'f+2 \delta f'\right]=0,
\end{split}
\end{equation*}
i.e., \Iff{}
\begin{align*}
 &\kappa f - 2 \delta g'= -\delta^{2} \kappa g+ \delta' f- 2\delta f'=0, \\
 &f \left(\kappa-\lambda\right)+\delta' g -2\delta g'=\delta^{2}g \left(\kappa-\lambda\right)-2 \delta'f+2 \delta f'=0.
\end{align*}
Treating the above system in the standard way we find that $$\lambda=\delta'=0.$$ If $$\kappa=0,$$ $\varPhi$ is right helicoid \footnote{The right helicoids are characterized by the conditions $\delta = c \!\in \! \mathbb{R}^{*}$ and $ \kappa  = \lambda  = 0$.}, $$f=c_{1}\in \mathbb{R^{*}}$$ and $$g=c_{2}\in \mathbb{R^{*}}.$$ If $$\kappa\neq0,$$ $\varPhi$ is orthoid of constant distribution parameter $c_{3}$, $$\kappa=\pm 2 g' \left| c_{4}-g^{2}\right|^{-1/2},\, c_{4}\in \mathbb{R^{*}},\, g'\neq0,\, g^{2}\neq c_{4}$$ and $$f = \pm c_{3} \left| c_{4}-g^{2}\right|^{1/2}.$$
So, we can state
 \begin{proposition}
   The support vector field $\overline{Q}$ of a right normalized skew ruled surface $\varPhi \subset \mathbb{E}^{3}$ is tangential to the curved asymptotic lines of $\varPhi$ \Iff{}\\
   \textup{(a)} $\varPhi$  is right helicoid, the function $f$ is a nonvanishing constant $c_{1}$ and the function $g$ is a nonvanishing constant $c_{2}$, or\\
   \textup{(b)} $\varPhi$  is orthoid of constant distribution parameter $c_{3}$ and the conical curvature and the function $f$ are given by
     \begin{equation*}
     \kappa=\pm 2 g' \left| c_{4}-g^{2}\right|^{-1/2},\, c_{4}\in \mathbb{R^{*}},\, g'\neq0,\, g^{2}\neq  c_{4}\,\, \textit{and}\,\, f = \pm c_{3} \left| c_{4}-g^{2}\right|^{1/2}.
     \end{equation*}
  \end{proposition}
From \eqref{555} and \eqref{572}, resp. \eqref{573}, we deduce: $\overline{Q}$ is tangential or orthogonal to the $\widetilde{K}$-curves \Iff{}
\begin{equation}\label{590}
  \begin{split}
      & 2\left(\kappa f -  \delta g'\right)v^{4}- \left(2\delta^{2}\kappa g - \delta' f + 2 \delta f'\right)v^{3}+\delta^{2}\left[2 f \left(\kappa-\lambda\right)+\delta'g-2\delta g'\right]v^{2} \\
       & -\delta^{2}\left[2 \delta^{2}g \left(\kappa-\lambda\right)-3 \delta'f+2 \delta f'\right]v-\delta^{4}\delta'g=0
  \end{split}
\end{equation}
or
\begin{equation}\label{595}
  \begin{split}
      & \delta' \left(\kappa f -  \delta g'\right)v^{3}+ \left[2 \delta^{2}f\left(\kappa \lambda + 1\right)-\delta^{2} \delta' \kappa g+\delta'^{2}f-\delta\delta'f'-2\delta^{3}\lambda g'\right]v^{2}\\
       & -\delta^{2}\left[2\delta^{2}g \left(\kappa\lambda+1\right)+\delta'\kappa f+2\lambda\left(\delta f'-\delta'f\right)-\delta \delta'g'\right]v\\
       &+\delta^{2}\delta'\left(\delta^{2}\kappa g-\delta'f+\delta f'\right)=0,
  \end{split}
\end{equation}
respectively. From \eqref{590} we have that $\overline{Q}$ is tangential to the $\widetilde{K}$-curves \Iff{}
\begin{equation*}
\begin{split}
    & \kappa f -  \delta g'=2\delta^{2}\kappa g - \delta' f + 2 \delta f'=2 f \left(\kappa-\lambda\right)+\delta'g-2\delta g'=0, \\
     & 2 \delta^{2}g \left(\kappa-\lambda\right)-3 \delta'f+2 \delta f'=\delta'=0,
\end{split}
\end{equation*}
from which we obtain that $\varPhi$ is relative minimal, $$\delta=c_{1}\in \mathbb{R^{*}},$$ $$\lambda=0$$ and $$f=\pm \left|c_{2}-c_{1}^{2}g^{2}\right|^{1/2},\, c_{2} \in \mathbb{R},\, c_{1}^{2}g^{2}\neq c_{2}.$$ From \eqref{595} we infer that $\overline{Q}$ is orthogonal to the $\widetilde{K}$-curves \Iff{}
\begin{equation*}
\begin{split}
    & \delta' \left(\kappa f -  \delta g'\right)=2 \delta^{2}f\left(\kappa \lambda + 1\right)-\delta^{2} \delta' \kappa g+\delta'^{2}f-\delta\delta'f'-2\delta^{3}\lambda g'=0, \\
     & 2\delta^{2}g \left(\kappa\lambda+1\right)+\delta'\kappa f+2\lambda\left(\delta f'-\delta'f\right)-\delta \delta'g'=\delta'\left(\delta^{2}\kappa g-\delta'f+\delta f'\right)=0,
\end{split}
\end{equation*}
that is, \Iff{} $\varPhi$ is relative minimal or $$\delta=c_{1}\in \mathbb{R^{*}}.$$ If $\varPhi$ is relative minimal we arrive at a contradiction.

If $$\delta=c_{1}\in \mathbb{R^{*}},$$ we have $$\kappa \lambda +1=\frac{c_{1} \lambda  g'}{f}$$ and $$f=\pm \left|c_{2}-c_{1}^{2}g^{2}\right|^{1/2},\, c_{2} \in \mathbb{R},\, c_{1}^{2}g^{2}\neq c_{2},$$ hence $$\kappa=\pm c_{1} g' \left|c_{2}-c_{1}^{2}g^{2}\right|^{-1/2}-\lambda^{-1},\, \lambda \neq 0.$$ Thus, we deduce
 \begin{proposition}
   The support vector field $\overline{Q}$ of a right normalized skew ruled surface $\varPhi \subset \mathbb{E}^{3}$ is\\
   \textup{(a)}tangential to the $\widetilde{K}$-curves of $\varPhi$ \Iff{} $\varPhi$ is an orthoid, relative minimal surface of constant distribution parameter $c_{1}$ and the function $f$ is given by
          \begin{equation*}
          f=\pm \left|c_{2}-c_{1}^{2}g^{2}\right|^{1/2},\, c_{2} \in \mathbb{R},\, c_{1}^{2}g^{2}\neq c_{2}.\\
          \end{equation*}
   \textup{(b)} orthogonal to the $\widetilde{K}$-curves of $\varPhi$ \Iff{} $\varPhi$  has constant distribution parameter $c_{1}$ and the conical curvature and the function $f$ are given by
          \begin{equation*}
           \kappa=\pm c_{1} g' \left|c_{2}-c_{1}^{2}g^{2}\right|^{-1/2}-\lambda^{-1},\,c_{2} \in \mathbb{R},\, \lambda \neq 0,\,c_{1}^{2}g^{2}\neq c_{2}\,\, \textit{and}\,\,f=\pm \left|c_{2}-c_{1}^{2}g^{2}\right|^{1/2}.
          \end{equation*}
  \end{proposition}
The following table summarizes the results:
\begin{footnotesize}
\begin{center}
  \begin{tabular}{|m{2.2cm}|m{4.5cm}|m{4.9cm}|}
    \hline
     \centering$\overline{Q}$ is \dots & \centering Type of the ruled surface $\varPhi$  & \centering \arraybackslash $f,\,g$ \\ [0.1cm]
    \hline \hline
    orthogonal to

    the generators & relative minimal & $f\!=\!\pm\delta\left|c\!-\!g^{2}\right|^{1/2},\, c\!\in \! \mathbb{R},\, g^{2}\!\neq\! c$ \\
    \hline
    tangential to

    the $u$-curves & orthoid, relative minimal  & $f\!=\!\pm \delta \left|c \!-\! g^{2}\right|^{1/2},\, c \!\in \! \mathbb{R}, \, g^{2}\!\neq\! c$ \\
    \hline
    orthogonal to

    the $u$-curves & $\kappa\!=\! \pm g' \left|c\!-\!g^{2}\right|^{-1/2}\!-\! \!\lambda^{-1}$,

    $c \!\in \! \mathbb{R}$, $\lambda\!\neq\!0$, $g^{2}\!\neq\! c$ & $f\!=\!\pm \delta \left|c\!-\!g^{2}\right|^{1/2}$\\
    \hline
    \multirow {2} {2cm}  {tangential to the curved asympt. lines} & right helicoid & $f\!=\!c_{1}\!\in \! \mathbb{R^{*}}$, $g\!=\!c_{2}\!\in \! \mathbb{R^{*}}$ \\
    \cline{2-3}
    & orthoid, $\delta \!=\! c_3 \in\mathbb{R^{*}}$,

     $\kappa\!=\!\pm 2 g' \left|c_{4}\!-\!g^{2}\right|^{-1/2}$,

     $c_{4}\!\in \! \mathbb{R^{*}}$, $g'\!\neq\! 0$, $g^{2}\!\neq\!  c_{4}$  & $f \!=\! \pm c_3 \left| c_{4}\!-\!g^{2}\right|^{1/2}$
   \\
    \hline
    tangential to

    the $\widetilde{K}$-curves & orthoid, relative minimal,

    $\delta \!=\! c_1 \in\mathbb{R^{*}}$ & $f\!=\!\pm \left|c_2\!-\!c_1^{2}g^{2}\right|^{1/2}$,

    $c_2 \!\in \! \mathbb{R}$, $c_1^2 g^{2}\!\neq\! c_2$ \\
    \hline
    orthogonal to

    the $\widetilde{K}$-curves & $\delta\!=\!c_{1}\!\in \! \mathbb{R^{*}}$,

    $\kappa\!=\!\pm c_{1} g' \left|c_2\! -\! c_{1}^{2}g^{2}\right|^{-1/2}\!-\!\!\lambda^{-1}$,

    $c_2 \!\in \! \mathbb{R}$, $\lambda \!\neq\! 0$, $c_{1}^{2}g^{2}\!\neq\! c_2$ & $f\!=\!\pm \left|c_2 \!-\! c_{1}^{2}g^{2}\right|^{1/2}$\\
    \hline
  \end{tabular}
\end{center}
\end{footnotesize}

\end{document}